\documentclass[11pt,a4paper]{article}

\usepackage{titlesec}
\usepackage{fancyhdr}
\usepackage{a4wide}
\usepackage{graphicx}
\usepackage{float}
\usepackage{amssymb}
\usepackage{amsmath}
\usepackage{amsfonts}
\usepackage{amsthm}
\usepackage{color}
\usepackage{mathrsfs}
\usepackage{array}
\usepackage{eucal}
\usepackage{tikz}
\usepackage[T1]{fontenc}
\usepackage{inputenc}
\usepackage[english]{babel}
\usepackage{lmodern}
\usepackage{hyperref}
\usepackage{geometry}
\usepackage{changepage}
\usepackage{bm}
\changepage{0pt}{}{}{}{}{0pt}{}{0pt}{6pt}
\usepackage[numbers]{natbib}
\usepackage{hyperref}
\hypersetup{
pdfpagemode=none,
pdftoolbar=true,        
pdfmenubar=true,        
pdffitwindow=false,     
pdfstartview={Fit},    
pdftitle={Abnormal acoustic transmission in a waveguide with perforated screens},    
pdfauthor={L. Chesnel, S.A. Nazarov},     
pdfsubject={},  
pdfcreator={L. Chesnel, S.A. Nazarov},   
pdfproducer={L. Chesnel, S.A. Nazarov}, 
pdfkeywords={}, 
pdfnewwindow=true,      
colorlinks=true,       
linkcolor=magenta,          
citecolor=red,        
filecolor=cyan,      
urlcolor=blue           
}

\newcommand{\dsp}{\displaystyle}

\newcommand{\eps}{\varepsilon}
\newcommand{\om}{\omega}
\newcommand{\Om}{\Omega}
\newcommand{\mrm}[1]{\mathrm{#1}}

\newcommand{\N}{\mathbb{N}}
\newcommand{\R}{\mathbb{R}}

\newcommand{\Capa}{\mrm{cap}}

\newtheorem{theorem}{Theorem}[section]

\newtheorem{proposition}[theorem]{Proposition}

\begin{document}

~\vspace{-0.5cm}
\begin{center}
{\sc \bf\huge Abnormal acoustic transmission\\[8pt] 
in a waveguide with perforated screens
}
\end{center}

\begin{center}
\textsc{Lucas Chesnel}$^1$, \textsc{Sergei A. Nazarov}$^{2}$\\[16pt]
\begin{minipage}{0.95\textwidth}
{\small
$^1$ INRIA/Centre de mathématiques appliquées, \'Ecole Polytechnique, Université Paris-Saclay, Route de Saclay, 91128 Palaiseau, France;\\
$^2$ St. Petersburg State University, Universitetskaya naberezhnaya, 7-9, 199034, St. Petersburg, Russia;\\[10pt]
E-mails: \texttt{lucas.chesnel@inria.fr}, \texttt{srgnazarov@yahoo.co.uk} \\[-14pt]
\begin{center}
(\today)
\end{center}
}
\end{minipage}
\end{center}
\vspace{0.4cm}

\noindent\textbf{Abstract.} We consider the propagation of the piston mode in an acoustic waveguide obstructed by two screens with small holes. In general, due to the features of the geometry, almost no energy of the incident wave is transmitted through the structure. The goal of this article is to show that tuning carefully the distance between the two screens, which form a resonator, one can get almost complete transmission. We obtain an explicit criterion, not so obvious to intuit, for this phenomenon to happen. Numerical experiments illustrate the analysis.\\[4pt]
\noindent\textbf{Key words.} Waveguides, perforated screens, asymptotic analysis, abnormal transmission.

\section{Introduction}

We study the propagation of acoustic waves in a 3D waveguide obtructed by two screens with holes of size $\eps$ where $\eps>0$ is a small parameter. We work at fixed frequency such that only the piston mode, constant in the transverse direction, can propagate. In this setting, the scattering of an incident piston mode is characterized by a reflection coefficient $R^{\eps}$ and a transmission coefficient $T^{\eps}$ (see \eqref{Field} below). Due to conservation of energy, we have 
\begin{equation}\label{ConservationNRJ}
|R^\eps|^2+|T^\eps|^2=1.
\end{equation}
In general, that is for arbitrary positions of the screens, due to the features of the geometry, almost no energy of the incident wave passes through the small holes and one observes almost complete reflection: $\lim_{\eps\to0}R^\eps=1$ and $\lim_{\eps\to0}T^\eps=0$. But tuning  carefully the distance $2L^{\eps}$ between the screens, we will see that one can get good energy transmission. More precisely, for certain choices of $L^0$, $L'$, $L''$ in $L^{\eps}:=L^0+\eps L'+\eps^2L''$ and for a certain condition \eqref{DefCriterion} on the shape of the holes, we can obtain $\lim_{\eps\to0}R^\eps=0$ and $\lim_{\eps\to0}T^\eps=T^0$ with $|T^0|=1$ (see Section \ref{sectionAnalysis}).

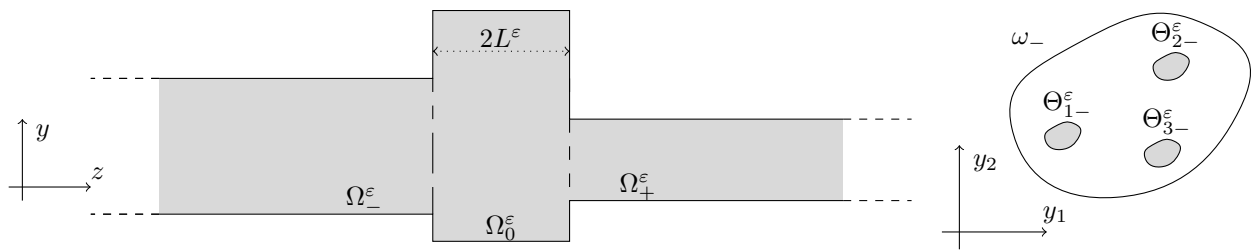
\begin{figure}[!ht]
\centering
\begin{tikzpicture}[scale=1.8]
\draw[fill=gray!30,draw=none](-2.5,0) rectangle (-0.5,1);
\draw[fill=gray!30,draw=none](2.5,0.1) rectangle (0.5,0.7);
\draw[fill=gray!30,draw=none](-0.5,-0.2) rectangle (0.5,1.5);
\draw[dashed] (-2.5,0)--(-3,0);
\draw[dashed] (-2.5,1)--(-3,1);
\draw[dashed] (2.5,0.1)--(3,0.1);
\draw[dashed] (2.5,0.7)--(3,0.7);
\draw[dashed] (-2.5,0)--(-3,0);
\draw[dashed] (-2.5,1)--(-3,1);
\draw (-2.5,0)--(-0.5,0)--(-0.5,-0.2)--(0.5,-0.2)--(0.5,0.1)--(2.5,0.1);
\draw (-2.5,1)--(-0.5,1)--(-0.5,1.5)--(0.5,1.5)--(0.5,0.7)--(2.5,0.7);
\draw (-0.5,0)--(-0.5,0.2);
\draw (-0.5,0.3)--(-0.5,0.6);
\draw (-0.5,0.7)--(-0.5,0.8);
\draw (-0.5,0.9)--(-0.5,1);
\draw (0.5,0.2)--(0.5,0.3);
\draw (0.5,0.4)--(0.5,0.5);
\draw (0.5,0.6)--(0.5,0.7);
\draw (0.5,0.9)--(0.5,1);
\node at (0,-0.1){\small $\Om_0^{\eps}$};
\node at (-1,0.1){\small $\Om_-^{\eps}$};
\node at (1,0.2){\small $\Om_+^{\eps}$};
\draw[dotted,<->] (-0.5,1.2)--(0.5,1.2);
\node at (0,1.3){\small $2L^{\eps}$};
\begin{scope}[shift={(-3.6,0)}]
\draw[->] (0,0.2)--(0.6,0.2);
\draw[->] (0.1,0.1)--(0.1,0.7);
\node at (0.65,0.3){\small $z$};
\node at (0.25,0.6){\small $y$};
\end{scope}
\end{tikzpicture}\quad
\begin{tikzpicture}[scale=2.3]
\begin{scope}[yshift=1cm]
\draw [fill=white] plot [smooth cycle, tension=1] coordinates {(-0.6,0.9) (0,0.5) (0.7,1) (0.5,1.5) (-0.2,1.4)};
\node at (-0.5,1.4){\small $\om_{-}$};
\begin{scope}[shift={(-0.32,0.7)},scale=0.15]
\draw [fill=gray!30] plot [smooth cycle, tension=1] coordinates {(-0.6,0.9) (0,0.5) (0.7,1) (0.5,1.5) (-0.2,1.4)};
\node at (0.25,2.2){\small $\Theta^{\eps}_{1-}$};
\end{scope}
\begin{scope}[shift={(0.3,1.1)},scale=0.15]
\draw [fill=gray!30] plot [smooth cycle, tension=1] coordinates {(-0.6,0.9) (0,0.5) (0.7,1) (0.5,1.5) (-0.2,1.4)};
\node at (0.25,2.2){\small $\Theta^{\eps}_{2-}$};
\end{scope}
\begin{scope}[shift={(0.25,0.6)},scale=0.15]
\draw [fill=gray!30] plot [smooth cycle, tension=1] coordinates {(-0.6,0.9) (0,0.5) (0.7,1) (0.5,1.5) (-0.2,1.4)};
\node at (0.25,2.2){\small $\Theta^{\eps}_{3-}$};
\end{scope}
\begin{scope}[shift={(-1,0.1)}]
\draw[->] (0,0.2)--(0.6,0.2);
\draw[->] (0.1,0.1)--(0.1,0.7);
\node at (0.65,0.3){\small $y_1$};
\node at (0.25,0.6){\small $y_2$};
\end{scope}
\end{scope}
\end{tikzpicture}
\caption{Side view of the waveguide $\Om^{\eps}$ (left) and picture of the left perforated screen (right). \label{Waveguide}} 
\end{figure}

\noindent Similar problems have been considered in \cite{BKNPS13,DeGr18} but with only one hole and with Dirichlet boundary conditions. In \cite{DeGr18}, the authors work with decompositions in Fourier series which are hard to generalize. In \cite{BKNPS13}, following \cite{MaNaPl}, the authors use techniques of construction of asymptotic expansions for boundary value problems in singularly perturbed domains. We also adopt this approach. But to the difference of the above mentioned publications, we employ a fine tuning procedure of the geometrical shape which allows us to reveal the complete transmission phenomenon.

\section{Setting}

First, we describe in detail the geometry (see Figure \ref{Waveguide}). Let $\om_{\pm}$, $\om_0$ be connected bounded open sets of $\R^2$ such that $\om_{\pm}\subset\om_0$. We define the 3D domains
\[
\Om^{\eps}_{\pm}:=\{x=(y,z)\in\R^3\,|\, y\in\om_{\pm}, \,\pm z > L^{\eps}\},\qquad \Om^{\eps}_{0}:=\om_0\times (-L^{\eps};L^{\eps}),
\]
with, for $\eps>0$ small,
\begin{equation}\label{DefL}
L^{\eps}=L^0+\eps L'+\eps^2L''.
\end{equation}
Here the parameters $L^0>0$, $L'$, $L''\in\R$ will be set later to observe interesting phenomena. Pick some points $\mathcal{P}_{j\pm}\in\om_{\pm}$, with $\mathcal{P}_{j\pm}\ne \mathcal{P}_{k\pm}$ for $j\ne k$, and some bounded sets $\theta_{j\pm}\subset\R^2$, $j=1,\dots, J_{\pm}$, with $J_{\pm}\in\N^{\ast}:=\{1,2,\dots\}$. Then define the small ``holes'' 
\[
\Theta_{j\pm}^{\eps}:=\{x=(y,z)\in\R^3\,|\, \eps^{-1}(y-\mathcal{P}_{j\pm})\in\theta_{j\pm},\,z=\pm L^{\eps}\}
\]
(see Figure \ref{Waveguide} right). Finally set
\[
\Om^{\eps}:=\Om^{\eps}_-\cup\Big(\bigcup_{j=1}^{J_{-}}\Theta_{j-}^{\eps}\Big)\cup \Om^{\eps}_0\cup\Big(\bigcup_{j=1}^{J_{+}}\Theta_{j+}^{\eps}\Big)\cup\Om^{\eps}_+.
\]
We consider the following problem with Neumann boundary condition
\begin{equation}\label{MainPb}
 \begin{array}{|rcll}
 \Delta u^\eps +\kappa^2  u^\eps&=&0&\mbox{ in }\Omega^\eps\\
 \partial_\nu u^\eps &=&0 &\mbox{ on }\partial\Omega^\eps.
\end{array}
\end{equation}
Here, $\Delta$ is the Laplace operator while $\partial_\nu$ corresponds to the derivative along the exterior normal. Furthermore, $u^\eps$ is the acoustic pressure in the medium while $\kappa>0$ is the wave number. Denote $\kappa_{\pm}^2$ the first positive eigenvalue of the Neumann Laplacian in $\om_{\pm}$. In (\ref{MainPb}), we work with $\kappa\in(0;\min(\kappa_-,\kappa_+))$ so that only the piston modes $\mrm{w}^{\mrm{in}}_{\pm}$, $\mrm{w}^{\mrm{out}}_{\pm}$ with
\begin{equation}\label{DefModes}
\mrm{w}^{\mrm{in}}_{\pm}(y,z)=\cfrac{e^{\mp i\kappa z}}{\sqrt{|\om_\pm|}}\,,\qquad \mrm{w}^{\mrm{out}}_{\pm}(y,z)=\cfrac{e^{\pm i\kappa z}}{\sqrt{|\om_\pm|}}\,,
\end{equation}
can propagate in $\Om^{\eps}_\pm$. Here $|\om_\pm|$ stands for the Lebesgue measure of the set $\om_{\pm}$. We are interested in the solution to the diffraction problem \eqref{MainPb} generated by the incoming wave $\mrm{w}^{\mrm{in}}_{-}$ in the trunk $\Om^{\eps}_-$. This solution admits the decomposition
\begin{equation}\label{Field}
 u^\eps(y,z)=\begin{array}{|ll}
\mrm{w}^{\mrm{in}}_{-}+R^\eps\,\mrm{w}^{\mrm{out}}_{-}(y,z+L^{\eps})+\dots \quad\mbox{ in }\Om^{\eps}_-\\[3pt]
\phantom{\mrm{w}^{\mrm{in}}_{-}+\ \ } T^\eps\,\mrm{w}^{\mrm{out}}_{+}(y,z-L^{\eps})+\dots \quad\mbox{ in }\Om_{+}^\eps
 \end{array}
\end{equation}
where $R^\eps\in\mathbb{C}$, $T^\eps\in\mathbb{C}$ are reflection and transmission coefficients. In this decomposition, the ellipsis stand for a remainder which decays at infinity with the rate $e^{-(\kappa_-^2-\kappa^2)^{1/2}|z|}$ in $\Om^{\eps}_-$ and $e^{-(\kappa_+^2-\kappa^2)^{1/2}|z|}$ in $\Om^\eps_+$. Note that the shifts $\pm L^{\eps}$ in the decomposition (\ref{Field}) are introduced to prepare the analysis below. With the normalisation (\ref{DefModes}), $R^{\eps}$ and $T^{\eps}$ satisfy the relation of conservation of energy (\ref{ConservationNRJ}). Our goal is to compute an asymptotic expansion of $R^{\eps}$, $T^{\eps}$ with respect to $\eps$ as $\eps$ tends to zero.

\section{Ansatz and auxiliary problems}
To observe interesting phenomena, we work with $L^0$ in (\ref{DefL}) such that
\[
L^0=\cfrac{\pi q}{2\kappa}\qquad\mbox{ where }q\in\N^{\ast}. 
\]
In this case, $\kappa^2$ is an eigenvalue of the problem
\begin{equation}\label{PbOm0}
 \begin{array}{|rcll}
-\Delta v &=&\lambda  v&\mbox{ in }\Omega^0_0:=\om_0\times(-L^0;L_0)\\
 \partial_\nu v &=&0 &\mbox{ on }\partial\Omega^0_0
\end{array}
\end{equation}
obtained by considering the limit $\eps\rightarrow0^+$ in the equation \eqref{MainPb} restricted to the resonator $\Om^{\eps}_0$. We shall assume that $\kappa^2$ is a simple eigenvalue and we denote by $\mrm{\bf v}$ the eigenfunction 
\begin{equation}\label{DefEigenFunc}
\mrm{\bf v}(x)=\cos(\kappa(z+L^0)).
\end{equation}
Far from the holes $\Theta^{\eps}_{j\pm}$, for the field $u^{\eps}$ in (\ref{Field}) we work with the ans\"atze
\begin{flalign}
u^\eps(x)&=\mrm{w}^{\mrm{in}}_{-}(y,z+L^{\eps})+R^0\,\mrm{w}^{\mrm{out}}_{-}(y,z+L^{\eps})+u^0_-(y,-(z+L^{\eps}))+\dots \quad\mbox{ in }\Om^{\eps}_- \label{AnsatzWaveguides1}\\[3pt]
u^\eps(x)&=T^0\,\mrm{w}^{\mrm{out}}_{+}(y,z-L^{\eps})+u^0_+(y,z-L^{\eps})+\dots \quad\mbox{ in }\Om^{\eps}_+ \label{AnsatzWaveguides2}\\[3pt]
u^\eps(x)&=\eps^{-1} a^0 \,\mrm{\bf v}(z)+\eps^0\,v'_{\eps}(x)+\eps v''(x)+\dots \quad\mbox{ in }\Om^{\eps}_0.\label{AnsatzWaveguides3}
\end{flalign}
Here $R^0$, $T^0$, $a^0$ are unknown complex constants and the functions  $u^0_{\pm}$, $v'_{\eps}$, $v''$ have to be determined. In particular, $u^0_{\pm}$ decay exponentially at infinity. The term $v'_{\eps}$ will depend on $\eps$ but this dependence will be rather explicit. In these expansions, the ellipsis stand for higher order terms which will be unimportant in the analysis.\\
\newline
In the vicinity of the holes $\Theta_{j\pm}^{\eps}$, we observe a boundary layer phenomenon. To capture it, we introduce the rapid variables $\xi_{j\pm}=(\xi_{j\pm}^1,\xi_{j\pm}^2,\xi_{j\pm}^3):=\eps^{-1}(x-\mathscr{P}^{\eps}_{j\pm})$ with $\mathscr{P}^{\eps}_{j\pm}:=(\mathcal{P}_{j\pm},\pm L^{\eps})$. We look for an expansion of $u^{\eps}$ in a neighbourhood of the holes $\Theta^{\eps}_{j\pm}$ of the form
\begin{equation}\label{ExpansionNearField}
u^{\eps}(x)=\eps^{-1}Z^{-1}_{j\pm}(\xi_{j\pm})+\eps^{0}Z^{0}_{j\pm}(\xi_{j\pm})+\dots,
\end{equation}
where the functions $Z^{-1}_{j\pm}$, $Z^{0}_{j\pm}$ are to determined. 
Observing that 
\[
(\Delta_x+\kappa^2)u^{\eps}(\eps^{-1}(x-\mathscr{P}^\eps_{j\pm}))=\eps^{-2}\Delta_{\xi_{j\pm}}u^{\eps}(\xi_{j\pm})+\dots,
\]
we are led to consider the Neumann problem 
\begin{equation}\label{PbBoundaryLayer}
-\Delta_\xi Z=0\quad\mbox{ in }\Xi_{j\pm},\quad\qquad \partial_\nu Z=0\quad\mbox{ on }\partial\Xi_{j\pm}
\end{equation}
where $\Xi_{j\pm}:=\R^3_-\cup\R^3_+\cup\theta_{j\pm}(0)$. Here, by convention, $\R^3_-:=\{\xi=(\xi_1,\xi_2,\xi_3)\in\R^2\times(-\infty;0)\}$, $\R^3_+:=\R^2\times(0;+\infty)$ and $\theta_{j\pm}(0):=\theta_{j\pm}\times\{0\}$. \\
\newline
Introduce $P_{j\pm}$ the capacity potential of the set $\theta_{j\pm}(0)$ which is defined as the solution to the problem
\[
-\Delta_{\xi}P_{j\pm}=0\quad\mbox{ in }\R^3\setminus\overline{\theta_{j\pm}(0)},\quad \qquad P_{j\pm}=1\quad\mbox{ on }\theta_{j\pm}(0),
\]
and decay at infinity. In the sequel, the asymptotic behaviour of $P_{j\pm}$ at infinity will play a major role. As $|\xi|\to+\infty$, we have (see e.g. \cite{Lank72})
\[
P_{j\pm}(\xi) = \frac{\Capa(\theta_{j\pm})}{|\xi|} + \vec{q}_{j\pm}\cdot\nabla\Phi(\xi) + O(|\xi|^{-3}),
\]
where $\Phi:=\xi\mapsto-1/(4\pi|\xi|)$ is the fundamental solution of the Laplace operator in $\R^3$ and $\vec{q}_{j\pm}$ is some given vector in $\R^3$. The term $\Capa(\theta_{j\pm})=(4\pi)^{-1}\int_{\R^3\setminus\theta_{j\pm}(0)}|\nabla P_{j\pm}|^2\,d\xi>0$ corresponds to the harmonic capacity \cite{PoSz51} of the planar crack $\theta_{j\pm}(0)$. Note that since $P_{j\pm}$ is even in $\xi_3$, we have $\vec{q}_{j\pm}=(\vec{q}^{\,1}_{j\pm},\vec{q}^{\,2}_{j\pm},0)$. Playing with symmetries, one can check that any smooth bounded solution of \eqref{PbBoundaryLayer} is of the form $c_0+c_1 W_{j\pm}(\xi)$ where $c_0$, $c_1$ are constants and where $W_{j\pm}$ is the function such that
\[
W_{j\pm}(\xi)=\left\{ \begin{array}{ll}
1-P_{j\pm}(\xi) &\xi_3>0 \\[3pt]
-1+P_{j\pm}(\xi) & \xi_3<0.
\end{array}
\right.
\]
Note that one can verify that $W_{j\pm}$ is harmonic and smooth in $\Xi_{j\pm}$. Moreover, $W_{j\pm}$ is odd in $\xi_3$. With this definition, for $\eta=\pm$, we have the expansion
\begin{equation}\label{DecompoW}
W_{j\eta}(\xi)=\pm 1 \pm 4\pi\Capa(\theta_{j\eta})\Phi(\xi) \mp \vec{q}_{j\eta}\cdot\nabla\Phi(\xi) + O(|\xi|^{-3}),\quad |\xi|\to+\infty,\,\pm \xi_3>0.
\end{equation}

\section{Asymptotic expansion of the scattering coefficients}

In order to identify the terms in the outer (\ref{AnsatzWaveguides1}), (\ref{AnsatzWaveguides2}), (\ref{AnsatzWaveguides3}) and inner (\ref{ExpansionNearField}) expansions of $u^{\eps}$, we will match the different behaviours in the neighbourhood of the holes $\Theta_{j\pm}^{\eps}$.\\
\newline 
$\star$ We start with the expansion (\ref{AnsatzWaveguides3}) of $u^{\eps}$ in $\Om^{\eps}_0$. From the Taylor formula and the expression (\ref{DefEigenFunc}) for $\mrm{\bf v}$, we have
\begin{equation}\label{Taylorv1}
\mrm{\bf v}(\pm L^{\eps}) = (\mp 1)^q(1-\eps^2(\kappa L')^2/2+O(\eps^3)).
\end{equation}
On the other hand, we observe that the expansions (\ref{AnsatzWaveguides1}), (\ref{AnsatzWaveguides2}) of $u^{\eps}$ in $\Om^{\eps}_{\pm}$ remain bounded as $\eps\to0$. Therefore, matching the constant behaviours at order $\eps^{-1}$, in the inner expansion (\ref{ExpansionNearField}), we get
\[
Z^{-1}_{j\pm}(\xi_{j\pm})=\frac{a^0}{2}\,(\mp 1)^{q}(1\mp  W_{j\pm}(\xi_{j\pm})).
\]
Note in particular that with this choice, $Z^{-1}_{j\pm}(\xi)$ indeed tends to zero as $|\xi|\to +\infty$, $\pm\xi_3>0$.\\
\newline
$\star$ Then we introduce the expansion (\ref{AnsatzWaveguides3}) of $u^{\eps}$ in $\Om^{\eps}_0$ in the initial problem and look at the terms of order $\eps^0$. This leads us to consider the problem
\begin{equation}\label{PbCorrection1}
\begin{array}{|lcl}
\dsp\Delta v'+\kappa^2v'=0\mbox{ in }\Om^{0}_0,\qquad \partial_\nu v'=0\mbox{ on }\partial\om_0\times(-L^0;L^0),\\[5pt]
\dsp\pm\partial_z v'(y,\pm L)= a^0(\mp 1)^{q}\kappa^2L'- a^0\pi(\mp 1)^{q}\sum_{j=1}^{J_{\pm}}\Capa(\theta_{j\pm})\delta(y-\mathcal{P}_{j\pm})\mbox{ for }y\in\om_{0}.
\end{array}
\end{equation}
To obtain the second boundary condition, we used the Taylor expansion
\begin{equation}\label{Taylorv2}
\partial_z\mrm{\bf v}(\pm L^{\eps}) = (\mp 1)^{1+q}\eps\kappa^2(L'+\eps L''+O(\eps^2)).
\end{equation}
It shows that the first term $\eps^{-1} a^0 \,\mrm{\bf v}(z)$ in (\ref{AnsatzWaveguides3}) generates an error of order $\eps^0$ on $\om_0\times\{\pm L^{\eps}\}$ which must be compensated. Moreover, the Dirac masses $\delta(y-\mathcal{P}_{j\pm})$ come from 
\begin{equation}\label{InnerFar}
Z^{-1}_{j\pm}(\xi)=\frac{a^0}{2}\,(\mp 1)^{q}(2+4\pi\Capa(\theta_{j\pm})\Phi(\xi)+\dots),\qquad |\xi|\to+\infty,\,\pm\xi_{3}>0.
\end{equation}
Since $\Phi(\xi)=-1/(4\pi|\xi|)$, note that $\eps^{-1}\Phi(\xi_{j\pm})=-1/(4\pi |x-\mathscr{P}^{\eps}_{j\pm}|)$ is a term of order $\eps^0$. We emphasize that $\mrm{\bf v}$ must have the singular behaviour of the Green's function at the points $(\mathcal{P}_{j\pm},\pm L^0)$. Multiplying the volume equation of (\ref{PbCorrection1}) by $\mrm{\bf v}$ and integrating twice by parts, we find that (\ref{PbCorrection1}) admits a solution if and only if $\int_{\partial \Om^0_0}\mrm{\bf v}\,\partial_{\nu}v'\,d\sigma=0$. For $a^0\ne0$, this is equivalent to have
\begin{equation}\label{PerturbL}
L'=\cfrac{\pi}{2\kappa^2|\om_0|}\,\sum_{\pm}\sum_{j=1}^{J_{\pm}}\Capa(\theta_{j\pm}).
\end{equation}
We emphasize that if $L'$ in (\ref{DefL}) is chosen different from the above value (\ref{PerturbL}), then we must have $a^0=0$. In this case, there is no term in $\eps^{-1}$ in (\ref{AnsatzWaveguides3}), (\ref{ExpansionNearField}) and we simply get almost complete reflection when $\eps$ tends to zero. Therefore, from now on, we assume that $L'$ is set as in (\ref{PerturbL}). Then the solution of (\ref{PbCorrection1}) is uniquely defined under the condition 
\[
\int_{\Om_0^0} v' \mrm{\bf v}\,dx=0.
\]
Since $L'>0$, we need to extend the function $v'$ defined in $\Om^0_0$ to $\Om^\eps_0$. We take $v_{\eps}'$ in the expansion (\ref{AnsatzWaveguides3}) by setting
\begin{equation}\label{ExpressionVPrime}
v_{\eps}'(y,z)=\begin{array}{|ll}
v'(y,z-\eps L'-\eps^2L'') & \mbox{ for }z>0\\
v'(y,z+\eps L'+\eps^2L'') & \mbox{ for }z<0.
\end{array}
\end{equation}
At $z=0$, $v_{\eps}'$ has the jumps
\begin{equation}\label{DefJumps}
\begin{array}{rcl}
[v_{\eps}'](y,0)&:=&v_{\eps}'(y,0^+)-v_{\eps}'(y,0^-)=-2\eps L'\partial_zv'(y,0)+O(\eps^2)\\[3pt]
\left[\partial_z v_{\eps}'\right](y,0)&:=&\partial_z v_{\eps}'(y,0^+)-\partial_z v_{\eps}'(y,0^-)=-2\eps L' \partial^2_{z}v'(y,0)+O(\eps^2).
\end{array}
\end{equation}
These jumps will be compensated with the term $v_{\eps}''$. The important point is that they occur in a region where $v'$ is smooth.\\
\newline
$\star$ The next step consists in matching the outer (\ref{AnsatzWaveguides1}), (\ref{AnsatzWaveguides2}) and inner (\ref{ExpansionNearField}) expansions of $u^{\eps}$ at order $\eps^0$ in $\Om^{\eps}_{\pm}$. In addition to (\ref{InnerFar}), we have
\[
Z^{-1}_{j\pm}(\xi)=-\frac{a^0}{2}\,(\mp 1)^{q}(4\pi\Capa(\theta_{j\pm})\Phi(\xi)+\dots),\qquad |\xi|\to+\infty,\,\mp\xi_3>0.
\]
As a consequence, we obtain that the functions 
$u^0_\pm$ in \eqref{AnsatzWaveguides1}, \eqref{AnsatzWaveguides2} must solve the following problems
\begin{equation}\label{DefOuterWaveguide}
\begin{array}{|l}
\Delta u^0_\pm+\kappa^2u^0_\pm=0\mbox{ in }\Om^{\sqsubset}_\pm:=\om_{\pm}\times(0;+\infty),\qquad \partial_{\nu} u^0_\pm=0\mbox{ on }\partial\om_{\pm}\times(0;+\infty)\\[3pt]
\dsp-\partial_zu^0_\pm(y,0)=i\kappa |\om_\pm|^{-1/2}S^0_{\pm}+a^0\pi\,(\mp 1)^{q}
\sum_{j=1}^{J_{\pm}}\Capa(\theta_{j\pm})\delta(y-\mathcal{P}_{j\pm})\mbox{ for }y\in\om_{\pm}.
\end{array}
\end{equation}
Here $S^0_{+}:=T^0$ and $S^0_{-}:=R^0-1$. With our choice for the ans\"atze, $u^0_\pm$ must be exponentially decaying at infinity. Multiplying (\ref{DefOuterWaveguide}) by $e^{i\kappa z}+e^{-i\kappa z}$ and integrating by parts, we have to impose that $\int_{\partial \Om^{\sqsubset}_\pm}(e^{i\kappa z}+e^{-i\kappa z})\,\partial_{\nu}u^0_\pm\,d\sigma=0$. This leads to the identities
\begin{equation}\label{relationRzeroTzero}
0=i\kappa |\om_\pm|^{+1/2} S^0_{\pm}+a^0\pi\,(\mp 1)^{q}\sum_{j=1}^{J_{\pm}}\Capa(\theta_{j\pm}).
\end{equation}
Introduce the generalized Green function $G_{j\pm}$ which solves
\begin{equation}\label{DefGreen}
\begin{array}{|l}
\Delta G_{j\pm}+\kappa^2G_{j\pm}=0\mbox{ in }\Om^{\sqsubset}_\pm,\qquad \partial_{\nu} G_{j\pm}=0\mbox{ on }\partial\om_{\pm}\times(0;+\infty)\\[3pt]
\dsp-\partial_zG_{j\pm}(y,0)=\delta(y-\mathcal{P}_{j\pm})-|\om_\pm|^{-1}\mbox{ for }y\in\om_{\pm}.
\end{array}
\end{equation}
Note that $G_{j\pm}$ is exponentially decaying at infinity (to show this, again multiply by $e^{i\kappa z}+e^{-i\kappa z}$ and integrate by parts). As $r_{j\pm}:=((y-\mathcal{P}_{j\pm})^2+z^2)^{1/2}$ tends to zero, we have the decomposition
\[
G_{j\pm}(x)=\cfrac{1}{2\pi r_{j\pm}}+\tilde{G}_{j\pm}(x)
\]
where the function $\tilde{G}_{j\pm}$ is smooth. The matrices $\mathcal{G}_{\pm}:=(\mathcal{G}^{\pm}_{jk})_{1 \le j,k\le J_{\pm}}$ with $\mathcal{G}^{\pm}_{jk}=\tilde{G}_{j\pm}(\mathcal{P}_{k\pm},0)$ is real and symmetric. With this notation, using identity (\ref{relationRzeroTzero}), we find that the functions $u^0_\pm$ introduced in (\ref{DefOuterWaveguide}) satisfy
\[
u^0_\pm=a^0\pi\,(\mp 1)^{q}\sum_{j=1}^{J_{\pm}}\Capa(\theta_{j\pm})G_{j\pm}.
\]
As a consequence, as $r_{j\pm}$ tends to zero, we have the representation
\begin{equation}\label{Devptu0}
u^0_\pm(x)=a^0\pi\,(\mp 1)^{q}\frac{\Capa(\theta_{j\pm})}{2\pi r_{j\pm}}+a^0\,(\mp 1)^{q}\,\mathcal{U}^0_{j\pm}+O(r_{j\pm})\quad\mbox{ with }\quad\mathcal{U}^0_{j\pm}:=\pi\,\sum_{k=1}^{J_{\pm}}\Capa(\theta_{k\pm})\mathcal{G}^{\pm}_{kj}.
\end{equation}
Now we define the terms $Z^{0}_{j\pm}$ in the near field expansions (\ref{ExpansionNearField}). From the expression (\ref{ExpressionVPrime}) of $v'_{\eps}$, as $r^{\eps}_{j\pm}:=|x-\mathscr{P}^{\eps}_{j\pm}|$ tends to zero, we obtain the expansion
\begin{equation}\label{Devptvp}
v'_{\eps}(x)=- a^0\pi(\mp 1)^{q}\,\frac{\Capa(\theta_{j\pm})}{2\pi r^{\eps}_{j\pm}}+a^0(\mp 1)^{q}\,\mathcal{V}'_{j\pm}+O(r^{\eps}_{j\pm})
\end{equation}
for some real constants $\mathcal{V}'_{j\pm}$ independent of $\eps$. 
Owing to (\ref{AnsatzWaveguides1}), (\ref{AnsatzWaveguides2}) and (\ref{Devptu0}), the function $Z^{0}_{j\pm}$ in (\ref{ExpansionNearField}) must verify
\[
Z^{0}_{j\pm}(\xi)=s^0_{\pm}+a^0\,(\mp 1)^{q}\,\mathcal{U}^0_{j\pm}+o(1),\qquad |\xi|\to+\infty,\,\pm\xi_3>0. 
\]
with $s^0_{+}:=T^0/|\om_+|^{1/2}$ and $s^0_{-}:=(1+R^0)/|\om_-|^{-1/2}$. Besides, owing to (\ref{AnsatzWaveguides3}) and (\ref{Devptvp}), we have
\[
Z^{0}_{j\pm}(\xi)=a^0(\mp 1)^{q}\,\mathcal{V}'_{j\pm}+o(1),\qquad |\xi|\to+\infty,\,\mp\xi_3>0. 
\]
We conclude that 
\[
Z^{0}_{j\pm}(\xi_{j\pm})=A_{j\pm} W_{j\pm}(\xi_{j\pm})+B_{j\pm}
\]
where, according to the decomposition (\ref{DecompoW}) of $ W_{j\pm}$, the constants $A_{j\pm}$, $B_{j\pm}$ solve the systems
\[
\pm A_{j\pm}+B_{j\pm}= s^0_{\pm}+a^0\,(\mp1)^{q}\,\mathcal{U}^0_{j\pm},\qquad \mp A_{j\pm}+B_{j\pm}= a^0\,(\mp1)^{q}\,\mathcal{V}'_{j\pm}
\]
Thus, we get
\[
A_{j\pm}=\pm(s^0_{\pm}+a^0\,(\mp1)^{q}\,\mathcal{U}^0_{j\pm}-a^0\,(\mp1)^{q}\,\mathcal{V}'_{j\pm})/2,\quad B_{j\pm}=(s^0_{\pm}+a^0\,(\mp1)^{q}\,\mathcal{U}^0_{j\pm}+a^0\,(\mp1)^{q}\,\mathcal{V}'_{j\pm})/2.
\]
$\star$ In (\ref{AnsatzWaveguides1})-(\ref{ExpansionNearField}), it only remains to define the term $v''$. Consider the problem
\begin{equation}\label{PbCorrection1Bis}
\begin{array}{|lcl}
\dsp\Delta v''+\kappa^2v''=0\mbox{ in }\Om^{0}_0\setminus(\om_0\times\{0\}),\qquad \partial_\nu v''=0\mbox{ on }\partial\om_0\times(-L^0;L^0),\\[5pt]
\left[v''\right](y,0)=2a^0L'\partial_zv'(y,0),\qquad\quad \left[\partial_z v''\right](y,0)=2 a^0L' \partial_{zz}v'(y,0)\\[4pt]
\dsp\pm\partial_z v''(y,\pm L)= a^0(\mp 1)^{q}\kappa^2L''\\[2pt]
\pm2\pi\dsp\sum_{j=1}^{J_{\pm}}A_{j\pm}\Capa(\theta_{j\pm})\delta(y-\mathcal{P}_{j\pm})\mp \frac{a^0}{2}(\mp 1)^q\dsp\sum_{j=1}^{J_{\pm}}\sum_{p=1,2}q_{j\pm}^p\cfrac{\partial\delta}{\partial y_p}(y-\mathcal{P}_{j\pm})\mbox{ for }y\in\om_{0}.
\end{array}
\end{equation}
Here the jumps at $z=0$ are introduced to compensate (\ref{DefJumps}). Moreover the boundary conditions of the third line have been obtained by using \eqref{Taylorv2} and by matching the expansions.\\
\newline 
Multiplying the volume equation of (\ref{PbCorrection1Bis}) by $\mrm{\bf v}$ and integrating twice by parts, we find that (\ref{PbCorrection1Bis}) admits a solution if and only if there holds a relation of the form (compatibility condition)
\begin{equation}\label{DefKappaAlpha}
(-1)^qT^0\,K_+
+(1+R^0)\,K_-+a_0(\alpha_1+\alpha_2L'')=0\qquad \mbox{ where }K_{\pm}:=\frac{\pi}{|\om_\pm|^{1/2}}\sum_{j=1}^{J_{\pm}}\Capa(\theta_{j\pm})
\end{equation}
and where $\alpha_1$, $\alpha_2$ are some real constants depending in particular on $\kappa$, $\om_0$, $\theta_{j\pm}$ but not on $\eps$ and $L''$. Thus together with (\ref{relationRzeroTzero}), we obtain the system
\begin{equation}\label{systemEquations}
\begin{array}{|lcl}
i\kappa \,T^0+a^0\,(- 1)^{q}K_+&=&0\\[3pt]
i\kappa \,(R^0-1)+a^0 K_-&=&0\\[3pt]
(-1)^qT^0\,K_+
+(1+R^0)\,K_-+a_0(\alpha_1+\alpha_2L'')&=&0.
\end{array}
\end{equation}
Solving \eqref{systemEquations}, we obtain the following proposition, the main result of this article.
\begin{proposition}\label{MainProposition}
Let $R^\eps$, $T^\eps$ be the scattering coefficients (see \eqref{Field}) in the geometry $\Om^\eps$ defined from the parameter $L^\eps=\pi q/(2\kappa)+\eps L'+\eps^2L''$. For $L'$ as in (\ref{PerturbL}), we have $\lim_{\eps\to0}R^\eps=R^0(L'')$, $\lim_{\eps\to0}T^\eps=T^0(L'')$ with 
\begin{equation}\label{ScaLim}
\dsp R^0(L'')=\cfrac{K_+^2-K_-^2-i\kappa\beta}{K_+^2+K_-^2-i\kappa\beta}\,,\qquad \mbox{and}\qquad T^0(L'')=\cfrac{2(-1)^{q+1}K_+K_-}{K_+^2+K_-^2-i\kappa\beta}\,.
\end{equation}
Here $\beta:=\alpha_1+\alpha_2 L''$ and $K_{\pm}$, $\alpha_1$, $\alpha_2$ are set in (\ref{DefKappaAlpha}). Besides, for the constant $a^0$ in (\ref{AnsatzWaveguides3}) we have
\begin{equation}\label{Formulaa0}
a^0=a^0(L'')=\cfrac{2i\kappa K_-}{K_+^2+K_-^2-i\kappa\beta}\,.
\end{equation}
\end{proposition}

\section{Analysis of the results}\label{sectionAnalysis}
\begin{figure}[!ht]
\centering
\raisebox{-0.2cm}{\begin{tikzpicture}[scale=0.95]
\draw[->] (-0.5,0) -- (5.1,0) node[right] {$\eps$};
\draw[->] (0,-0.2) -- (0,3.3) node[left] {$L^\eps$};
\node at (-0.1,3.14/2-0.3) [left] {\small $\cfrac{\pi q}{2\kappa}$};
\begin{scope}
\clip(-2,-0.5) rectangle (4.5,3.2);
\draw[domain=0.01:4.7,smooth,variable=\x,green!60!black] plot ({\x},{3.14/2+\x-0.5*\x*\x});
\draw[domain=0.01:4.7,smooth,variable=\x,green!60!black] plot ({\x},{3.14/2+\x-0.4*\x*\x});
\draw[domain=0.01:4.7,smooth,variable=\x,green!60!black] plot ({\x},{3.14/2+\x-0.3*\x*\x});
\draw[domain=0.01:4.7,smooth,variable=\x,green!60!black] plot ({\x},{3.14/2+\x-0.2*\x*\x});
\draw[domain=0.01:4.7,smooth,variable=\x,green!60!black] plot ({\x},{3.14/2+\x-0.1*\x*\x});
\draw[domain=0.01:4.7,smooth,variable=\x,green!60!black] plot ({\x},{3.14/2+\x+0*\x*\x});
\draw[domain=0.01:4.7,smooth,variable=\x,green!60!black] plot ({\x},{3.14/2+\x+0.1*\x*\x});
\draw[domain=0.01:4.7,smooth,variable=\x,green!60!black] plot ({\x},{3.14/2+\x+0.2*\x*\x});
\draw[domain=0.01:4.7,smooth,variable=\x,green!60!black] plot ({\x},{3.14/2+\x+0.3*\x*\x});
\draw[domain=0.01:4.7,smooth,variable=\x,green!60!black] plot ({\x},{3.14/2+\x+0.4*\x*\x});
\draw[domain=0.01:4.7,smooth,variable=\x,green!60!black] plot ({\x},{3.14/2+\x+0.5*\x*\x});
\draw[domain=0.01:4.7,smooth,variable=\x,green!60!black] plot ({\x},{3.14/2+\x+0.6*\x*\x});
\draw[domain=0.01:4.7,smooth,variable=\x,green!60!black] plot ({\x},{3.14/2+\x+0.7*\x*\x});
\end{scope}
\draw[magenta,dashed,very thick] (0.4,-0.2) -- (0.4,3.3);
\node at (0.4,-0.2) [below] {$\eps_0$};
\end{tikzpicture}}\qquad\qquad\qquad
\begin{tikzpicture}[scale=1.8]
\draw[thick] (0,0) circle (1);
\draw[red, thick] (-0.5,0) circle (0.5);
\draw[blue, thick] (0.5,0) circle (0.5);
\draw[->] (0,-1.05)--(0,1.1);
\draw[->] (-1.05,0)--(1.1,0);
\fill[red] (-0.5,0.5) circle (0.4mm);
\node at(0.5,0.5){\textcolor{blue}{{\scriptsize $\square$}}};
\end{tikzpicture}
\caption{Left: paths $\{\gamma_{L''}(\eps)=(\eps,\pi q/(2\kappa)+\eps L'+\eps^2L'',\,\eps>0\}\subset\mathbb{R}^2$ for several values of $L''$.  According to the chosen path, the limit of the scattering coefficients along this path as $\eps\to0^+$ is different. With this picture, we understand why for a fixed small $\eps_0$, the scattering coefficients have a rapid variation as the distance between the screens changes in a vicinity of $\pi q/(2\kappa)$. Right: sets $\{R^0(L'')\,,\,L''\in\mathbb{R}\}$ (\textcolor{blue}{{\scriptsize $\square$}}) and $\{T^0(L'')\,,\,L''\in\mathbb{R}\}$ (\raisebox{0.5mm}{\protect\tikz \fill[red] (0,3) circle (0.6mm);}) in the complex plane where $R^0(L'')$, $T^0(L'')$ are defined in (\ref{ScaLim}). Here $q$ is odd and $K_-=K_+$. The black bold line represents the unit circle. \label{Illustration}}
\end{figure}
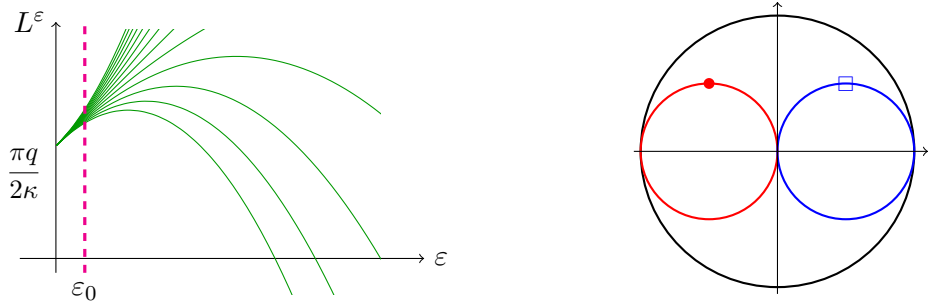

First, we observe that the coefficients $R^0(L'')$, $T^0(L'')$ defined in Proposition \ref{MainProposition} satisfy the relation of conservation of energy $|R^0(L'')|^2+|T^0(L'')|^2=1$. On the other hand, $R^0$ vanishes for a certain $L''=L''_{\diamond}$ (such that $\beta(L''_{\diamond})=0\Leftrightarrow \alpha_1+\alpha_2L''_{\diamond}=0$) if and only if $K_-=K_+$. This is equivalent to
\begin{equation}\label{DefCriterion}
\cfrac{1}{|\om_-|^{1/2}}\,\sum_{j=1}^{J_{-}}\Capa(\theta_{j-})=\cfrac{1}{|\om_+|^{1/2}}\,\sum_{j=1}^{J_{+}}\Capa(\theta_{j+}).
\end{equation}
Then we have $T^0(L''_{\diamond})=(-1)^{q+1}$ and $a^0(L''_{\diamond})=i\kappa/K_-$. Note that in order (\ref{DefCriterion}) to be satisfied, we do not need $J_-=J_+$ or $\om_-=\om_+$. Moreover, the position of the holes does not play any role because we deal with the piston modes. When (\ref{DefCriterion}) is met, setting $\tilde{\beta}=\kappa\beta/(2K^2_-)$, we obtain
\begin{equation}\label{resultsCircle}
\dsp R^0(L'')=\cfrac{-i\tilde{\beta}}{1-i\tilde{\beta}}\,,\qquad\qquad T^0(L'')=\cfrac{(-1)^{q+1}}{1-i\tilde{\beta}}\,,\qquad\qquad a^0(L'')=\cfrac{i\kappa /K_-}{1-i\tilde{\beta}}\,.
\end{equation}
In this case, there holds $R^0(L'')+(-1)^{q+1}T^0(L'')=1$ and as $L''$ varies in $\R$, $R^0(L'')$ runs on the circle centred at $1/2$ of radius $1/2$ while $T^0$ runs on the circle centred at $(-1)^{q+1}/2$ of radius $1/2$ (see Figure \ref{Illustration} right).

\noindent When the geometry is symmetric with respect to the plane $z=0$, we have $K_-=K_+$ and so $R^0(L''_{\diamond})=0$, $T^0(L''_{\diamond})=(-1)^{q+1}$. But in this situation, working with symmetries for example as in \cite{ChNa18}, we can get better and show that for $\eps>0$ small enough, there is $L^{\eps}$ close to $\pi q/(2\kappa)+\eps L'+\eps^2L''_{\diamond}$ such that $R^\eps=0$ and $T^\eps=(-1)^{q+1}$ (exactly and not asymptotically). We stress that for exact complete transmission, the position of holes (and not only their shapes and numbers) matters.

\section{Numerical illustrations}
\begin{figure}[!ht]
\centering
\raisebox{0cm}{\includegraphics[width=0.47\textwidth]{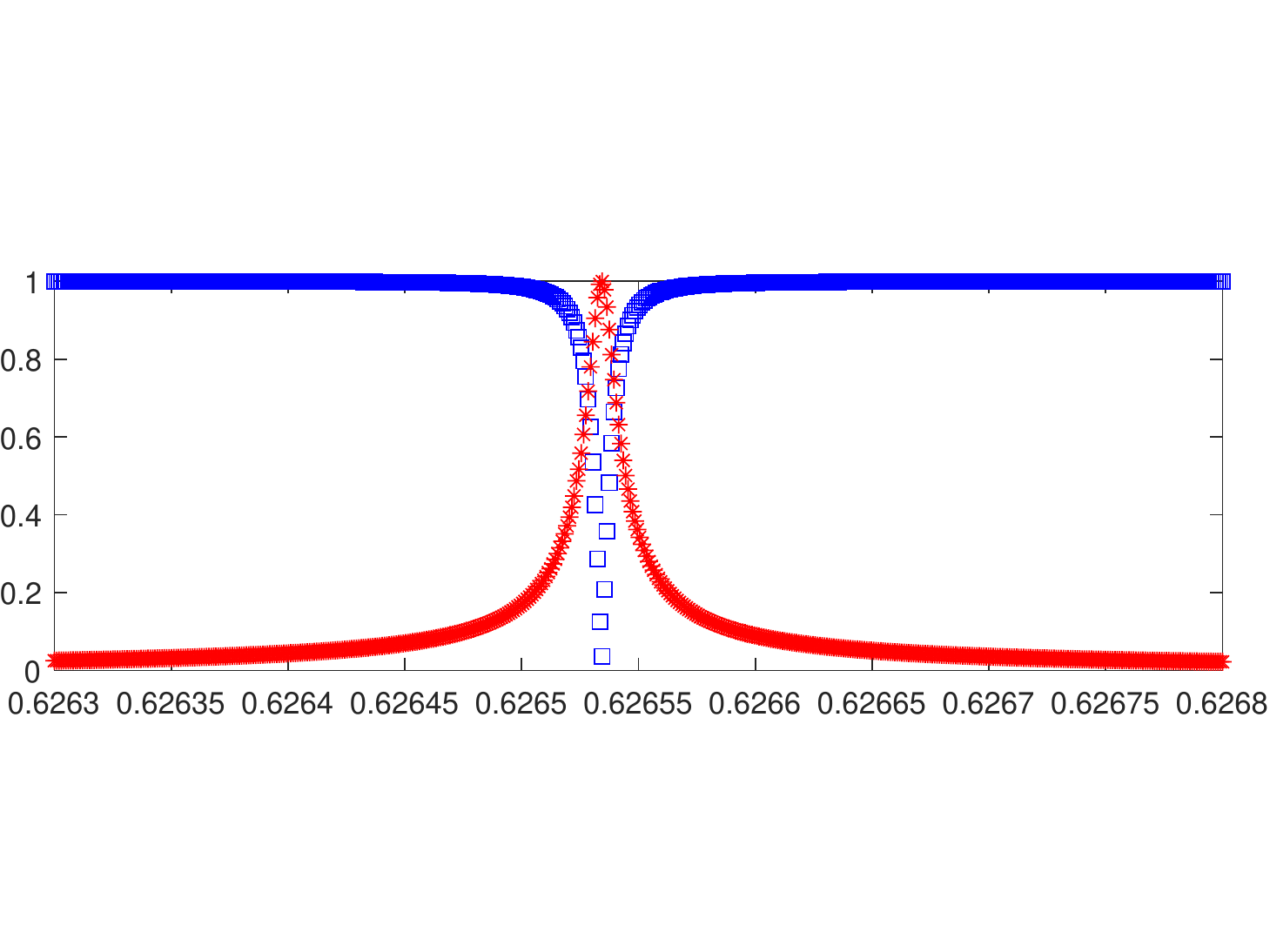}}\qquad\qquad
\includegraphics[trim={1.4cm 0cm 1.8cm 0cm},clip,width=0.4\textwidth]{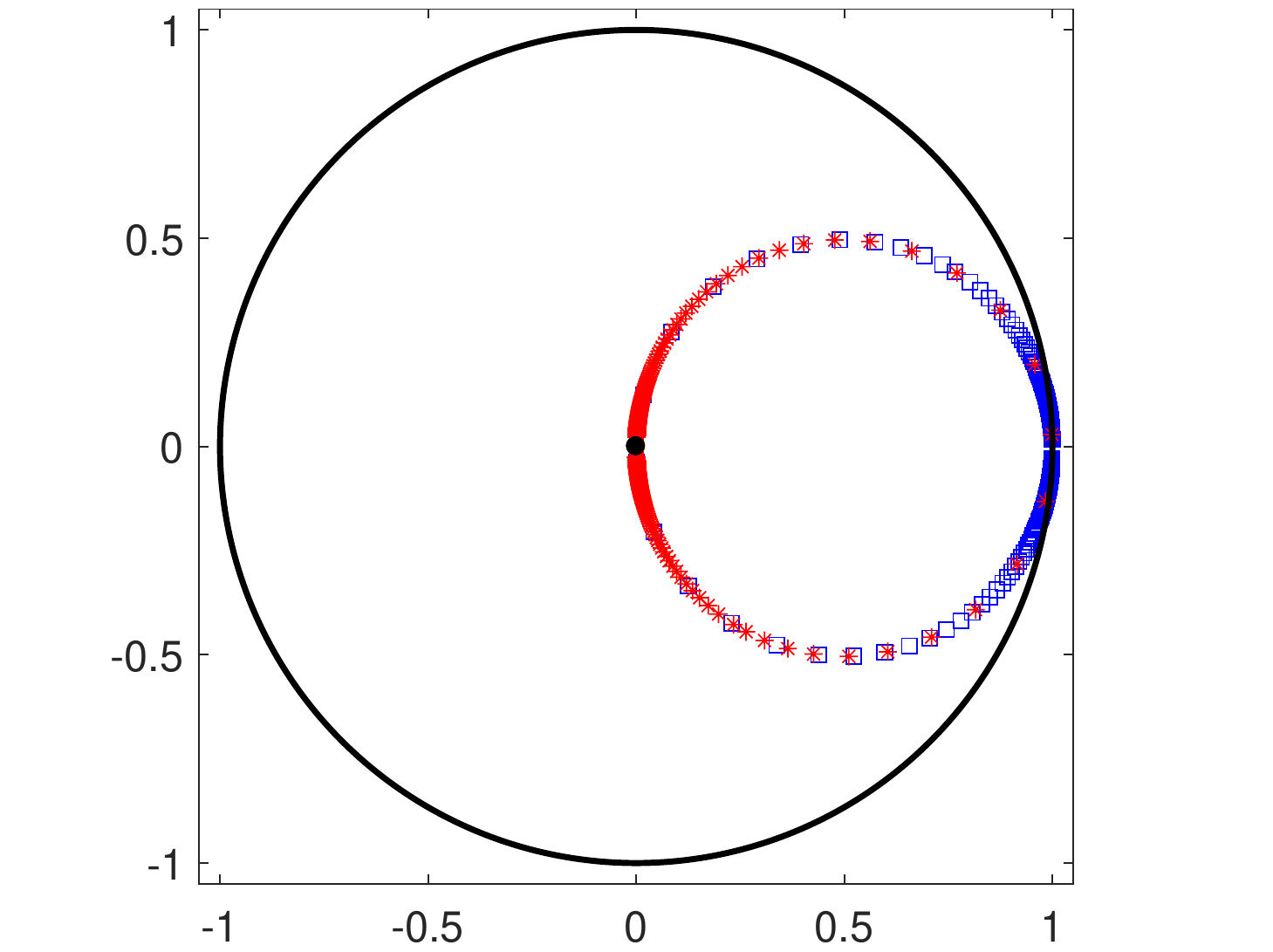}
\caption{Left: curves $L\mapsto |R^{\eps}(L)|$ (\textcolor{blue}{{\scriptsize $\square$}}) and $L\mapsto |T^{\eps}(L)|$ (\raisebox{0.3mm}{\textcolor{red}{\scriptsize{$\rlap{+}{\times}$}}}). Right: $L\mapsto R^{\eps}(L)$ (\textcolor{blue}{{\scriptsize $\square$}}) and $L\mapsto T^{\eps}(L)$ (\raisebox{0.3mm}{\textcolor{red}{\scriptsize{$\rlap{+}{\times}$}}}) in the complex plane. According to the conservation of energy, we have $|R^{\eps}(L)|^2+|T^{\eps}(L)|^2=1$. Therefore the scattering coefficients are located inside the unit disk delimited by the black bold line. For both pictures, $L$ takes values close to $L^0=\pi/(2\kappa)=0.625$ and $\eps=10^{-4}$. 
\label{MatriceScattering_sym}}
\end{figure}

\begin{figure}[!ht]
\centering
\raisebox{0.0cm}{\includegraphics[width=0.47\textwidth]{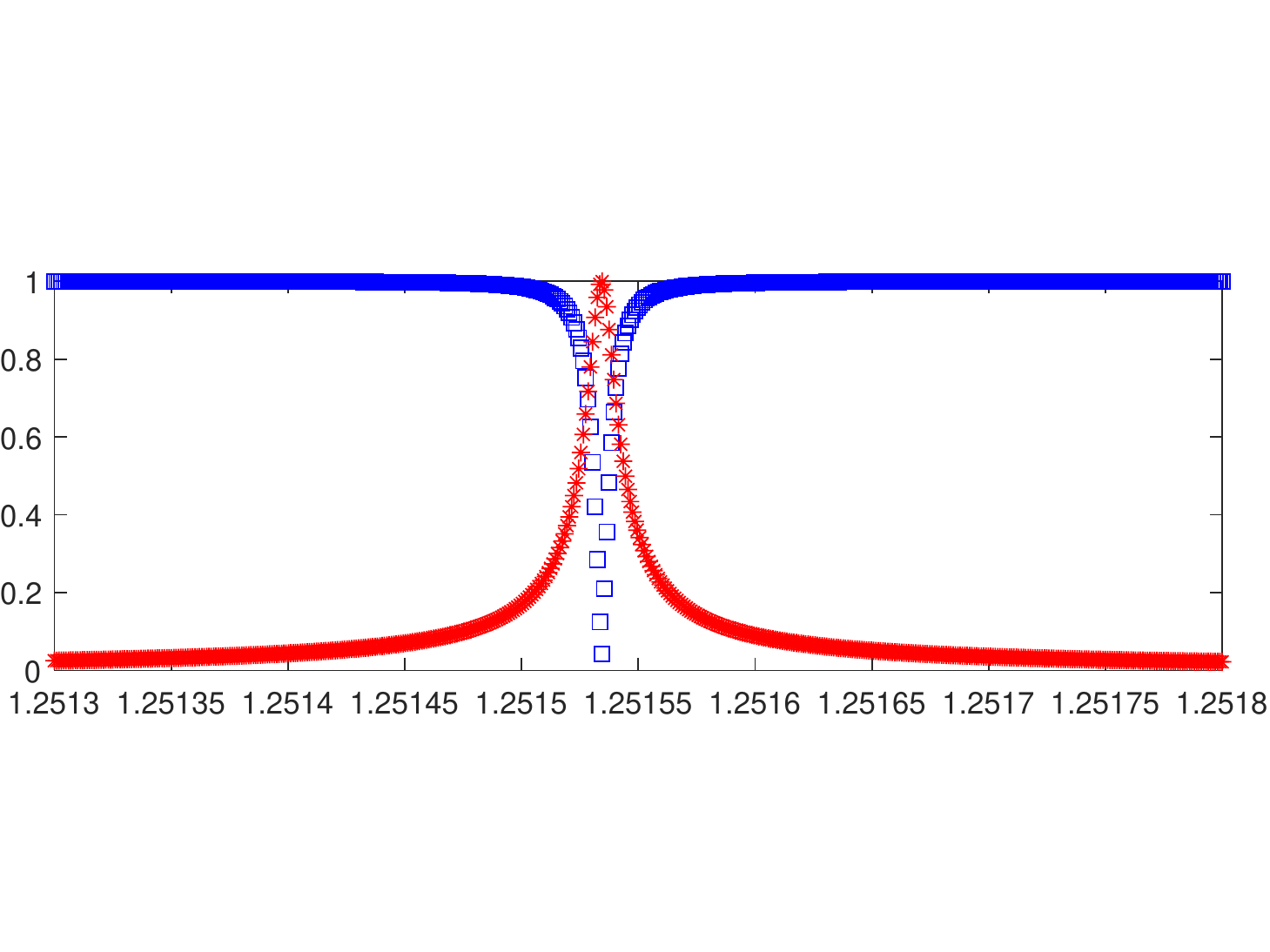}}\qquad\quad
\includegraphics[trim={1.4cm 0cm 1.8cm 0cm},clip,width=0.4\textwidth]{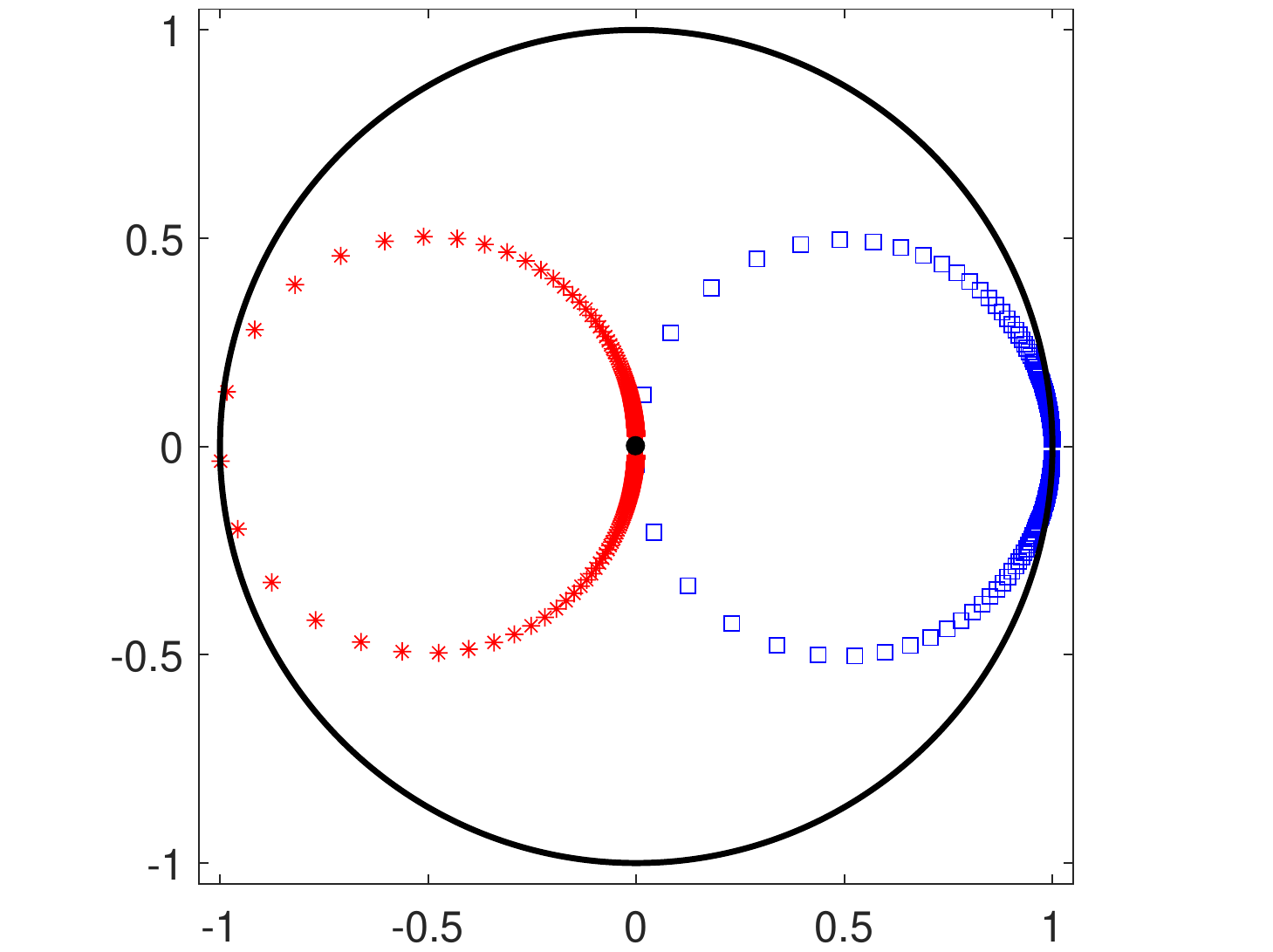}
\caption{Same quantities as in Figure \ref{MatriceScattering_sym} but with $L$ taking values close to $L^0=2\pi/(2\kappa)=1.25$. 
\label{MatriceScattering_symq2}}
\end{figure}

\begin{figure}[!ht]
\centering
\begin{tabular}{l}
\raisebox{6.5mm}{1)}\ \includegraphics[width=0.793\textwidth]{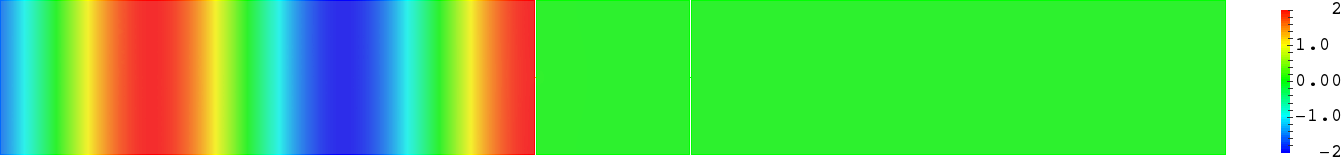}\\
\raisebox{6.5mm}{2)}\ \includegraphics[width=0.8\textwidth]{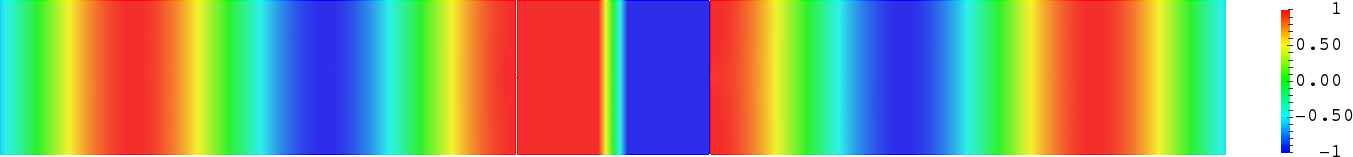}\\
\raisebox{6.5mm}{3)}\ \includegraphics[width=0.8\textwidth]{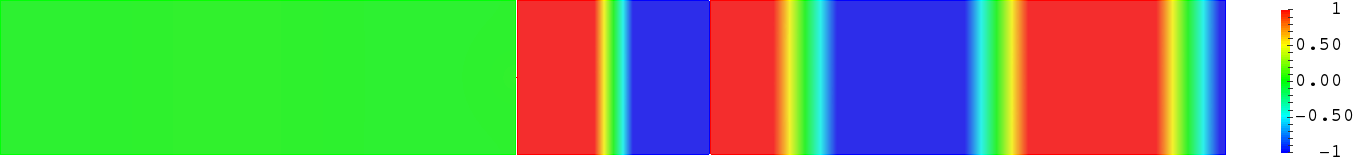}\\
\raisebox{6.5mm}{4)}\ \includegraphics[width=0.8\textwidth]{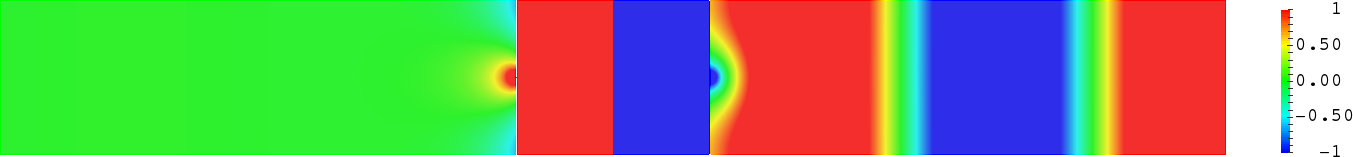}
\end{tabular}
\caption{1) $\Re e\,u^{\eps}$ for $L=0.5$. 2) $\Re e\,u^{\eps}$ for $L=0.6265\approx L^{\star}$. 3) (resp. 4)) $\Re e\,(u^{\eps}-u_i)$ (resp. $\Im m\,(u^{\eps}-u_i)$) (scattered field) for $L=0.6265\approx L^{\star}$. Here $u_i(x)=\mrm{w}^{\mrm{in}}_{-}(z+L)$. 
\label{Fields}}
\end{figure}

In this section, we illustrate the results we have obtained above. To simplify the numerical  implementation, we work in 2D. We emphasize that the asymptotic analysis is different from the above 3D setting. However, the physical phenomena are similar. For the experiments, we define  the waveguide $\Omega^{\eps}$ such that for $L>0$, 
\begin{equation}\label{DefGuideNum}
\Omega^{\eps}=\R\times(0;1)\setminus \{\Sigma_-^{\eps}\cup\Sigma_+^{\eps}\}\qquad\mbox{ with }\qquad \Sigma_{\pm}^{\eps}:=\{\pm L\}\times I_\pm^{\eps}.
\end{equation}
Here the sets $I_{\pm}^{\eps}$ depend on the situation and will be given below. We take $\kappa=0.8\pi<\kappa_{\pm}=\pi$ so that only the piston modes (see (\ref{DefModes})) can propagate. We compute numerically the scattering solution $u^{\eps}$ defined in (\ref{Field}). To proceed, we use a $\mrm{P2}$ finite element method in a truncated geometry. On the artificial boundary created by the truncation, a Dirichlet-to-Neumann operator with 15 terms serves as a transparent condition. Once we have computed $u^{\eps}$, it is easy to obtain the scattering coefficients $R^{\eps}$, $T^{\eps}$  in the representation (\ref{Field}). For the numerics, we take $\eps=10^{-4}$.\\
\newline
For the numerics of Figure \ref{MatriceScattering_sym}, in (\ref{DefGuideNum}) we take $I_-^{\eps}=I_+^{\eps}=(0;1)\setminus[1/2-\eps/2;1/2+\eps/2]$ (the holes are centered on the middle line of the waveguide). For $\kappa=0.8\pi$, the first critical length is  $L^0=\pi/(2\kappa)=1/1.6=0.625$. In Figure \ref{MatriceScattering_sym}, we display the scattering coefficients for $L$ varying close to $0.625$. 
As expected, when $\eps$ is small, for most values of $L$, the energy of the incident field is almost completely backscattered and the transmission coefficient $T$ is close to zero. In accordance with the discussion of Section \ref{sectionAnalysis} (remark that $K_-=K_+$ and even strongly, the geometry is symmetric with respect to $z=0$), we observe the phenomenon of complete transmission for some $L=L^{\star}$. As expected (see formula (\ref{PerturbL})), we note that $L^{\star}>L^0$. In Figure \ref{MatriceScattering_sym}, we find back the circles characterised by the formulas (\ref{resultsCircle}) for the asymptotic behaviour of the scattering coefficients. In Figure \ref{MatriceScattering_symq2}, we display the same quantities as in Figure \ref{MatriceScattering_sym} but with $L$ varying close to the second critical length  $L^0=2\pi/(2\kappa)=1.25$. Again, we get results in agreement with (\ref{resultsCircle}).\\
\newline
In Figure \ref{Fields}, we display the field $u^{\eps}$ for two different values of $L$, namely for a generic one where $T$ is almost zero and for $L\approx L^{\star}$. For $L\approx L^{\star}$, we indeed observe that the scattering field is exponentially decaying in the incident branch. For $L\approx L^{\star}$, we also note that the imaginary part of $u^{\eps}$ is large in the resonator, of the order $\eps^{-1}$. This is coherent with the formula \eqref{Formulaa0} which indicates that $a^0(L^{\star})$ is purely imaginary.

\begin{figure}[!ht]
\centering
\includegraphics[trim={1.4cm 0cm 1.8cm 0cm},clip,width=0.4\textwidth]{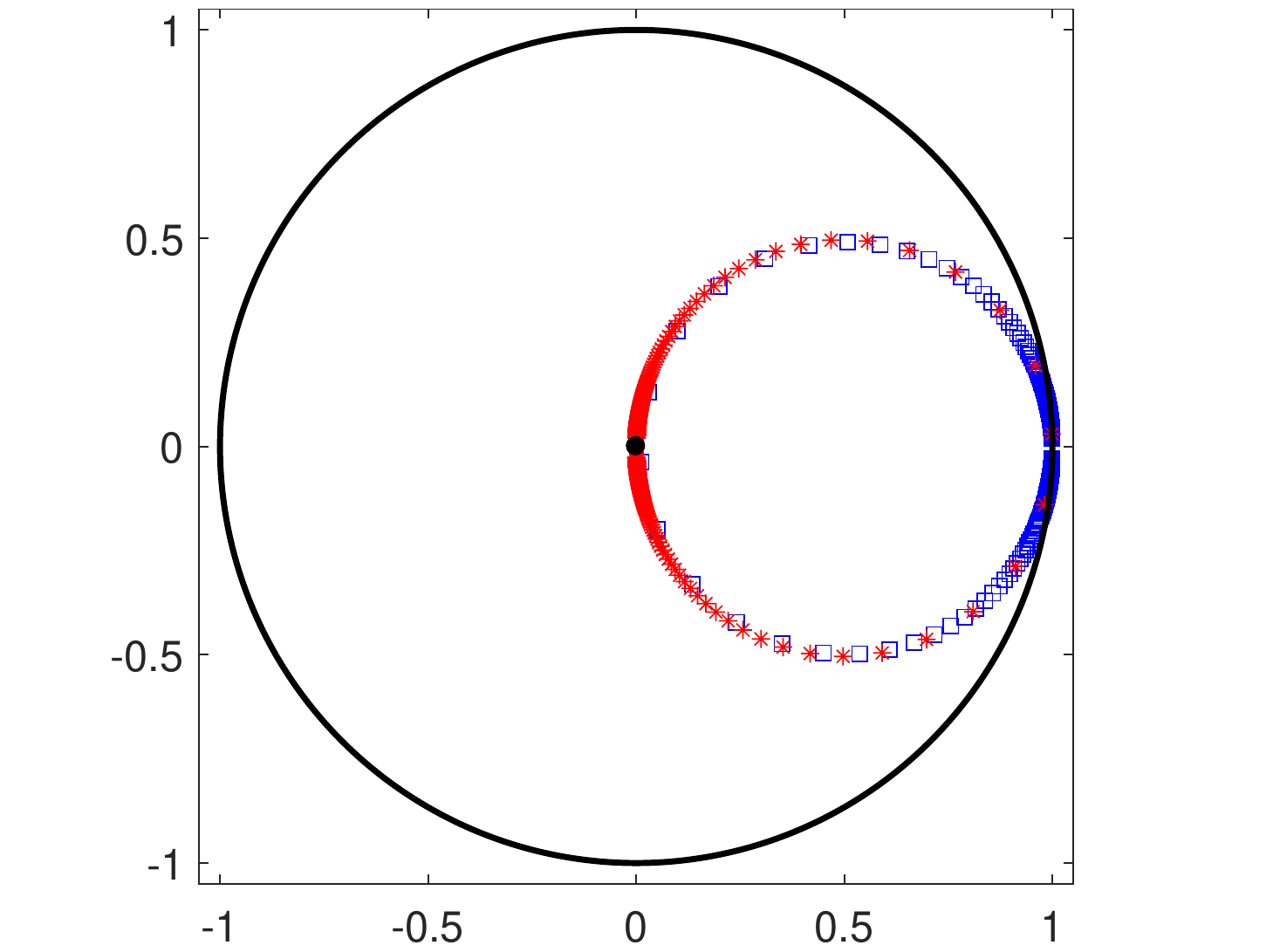}\qquad
\includegraphics[trim={1.4cm 0cm 1.8cm 0cm},clip,width=0.4\textwidth]{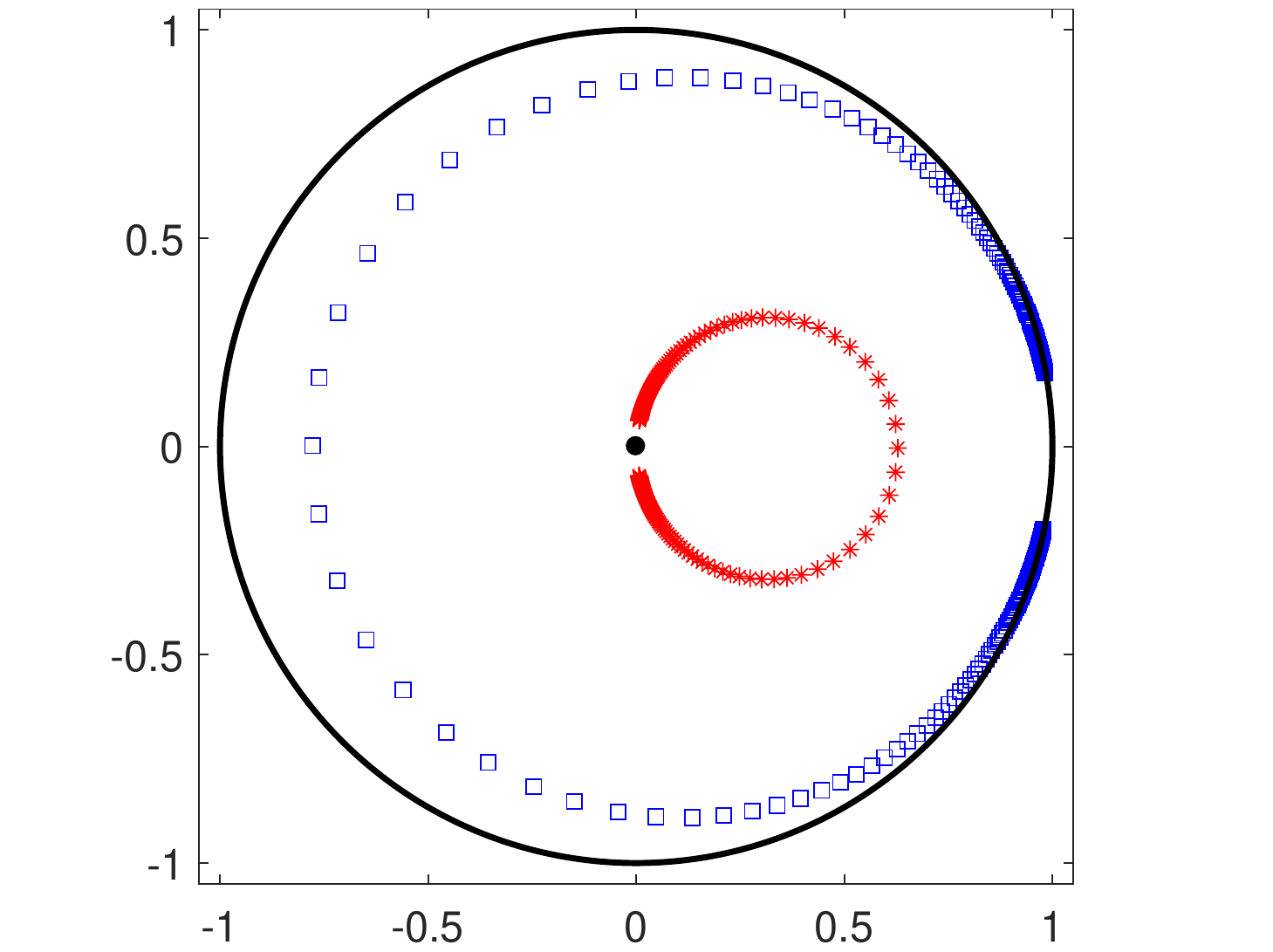}
\caption{Same quantities as in Figure \ref{MatriceScattering_sym} but in the geometries defined by (\ref{defGeomDeca}) (left) and (\ref{defGeomLarge}) (right).
\label{MatriceScattering_decalage}}
\end{figure}

\newpage

\noindent  For the numerics of Figure \ref{MatriceScattering_decalage} left, in (\ref{DefGuideNum}) we take 
\begin{equation}\label{defGeomDeca}
I_-^{\eps}=(0;1)\setminus[0.1-\eps/2;0.1+\eps/2]\qquad\mbox{ and }\qquad I_+^{\eps}=(0;1)\setminus[0.7-\eps/2;0.7+\eps/2].
\end{equation}
In this case, the holes are not at the center of the waveguide and there is no symmetry with respect to $z=0$. However, we still have $K_+=K_-$ and (\ref{resultsCircle}) indicates that we should observe almost complete transmission for a certain $L^{\ast}$. And this is what we get. For the numerics of Figure \ref{MatriceScattering_decalage} right, we take
\begin{equation}\label{defGeomLarge}
I_-^{\eps}=(0;1)\setminus[0.5-3\eps/2;0.5+3\eps/2]\qquad\mbox{ and }\qquad I_+^{\eps}=(0;1)\setminus[0.5-\eps/2;0.5+\eps/2].
\end{equation}
In other words, the right hole is three times larger than the left hole. In this case, we have $\Capa(\theta_{j-})=3\Capa(\theta_{j+})$ so that $K_-=3K_+$ and $\inf_{L}|R^0(L)|=4/5$ according to (\ref{resultsCircle}). This is indeed what we observe.

\section*{Acknowledgments} 
The research of S.A. Nazarov was supported by the grant No. 17-11-01003 of the Russian Science Foundation.

\bibliography{Bibli}

\def\cprime{$'$}
\begin{thebibliography}{1}

\bibitem{BKNPS13}
{L.M.} Baskin, M.~Kabardov, P.~Neittaanm{\"a}ki, {B.A.} Plamenevskii, and
  {O.V.} Sarafanov.
\newblock Asymptotic and numerical study of resonant tunneling in
  two-dimensional quantum waveguides of variable cross section.
\newblock {\em Comput. Math. Math. Phys.}, 53(11):1664--1683, 2013.

\bibitem{ChNa18}
L.~Chesnel and {S.A.} Nazarov.
\newblock Non reflection and perfect reflection via {Fano} resonance in
  waveguides.
\newblock {\em Comm. Math. Sci.}, 16(7):1779--1800, 2018.

\bibitem{DeGr18}
A.~Delitsyn and {D.S.} Grebenkov.
\newblock Mode matching methods for spectral and scattering problems.
\newblock {\em The Quarterly Journal of Mechanics and Applied Mathematics},
  71(4):537--580, 2018.

\bibitem{Lank72}
N.~S. Landkof.
\newblock {\em Foundations of modern potential theory}.
\newblock Springer-Verlag, New York-Heidelberg, 1972.
\newblock Translated from the Russian by A. P. Doohovskoy, Die Grundlehren der
  mathematischen Wissenschaften, Band 180.

\bibitem{MaNaPl}
{V.G.} {Maz'ya}, {S.A.} Nazarov, and {B.A.} Plamenevski{\u\i}.
\newblock {\em {Asymptotic theory of elliptic boundary value problems in
  singularly perturbed domains, Vol. 1 \ 2}}.
\newblock {Birkh\"{a}user}, Basel, 2000.
\newblock Translated from the original German 1991 edition.

\bibitem{PoSz51}
G.~P{\'o}lya and G.~Szeg{\"o}.
\newblock {\em Isoperimetric {I}nequalities in {M}athematical {P}hysics}.
\newblock Annals of Mathematics Studies, no. 27. Princeton University Press,
  Princeton, N. J., 1951.

\end{thebibliography}
\bibliographystyle{plain}
\end{document}